\documentclass[12pt,reqno]{amsart}
\usepackage[cp1251]{inputenc}
\usepackage[english]{babel}
\usepackage{amsfonts}
\usepackage{amssymb,amsmath,amsthm,eufrak,euscript}
\usepackage{graphicx}
\usepackage[usenames]{color}
\usepackage{colortbl}
\newtheorem{lemma}{Lemma}
\newtheorem{theorem}{Theorem}

\newtheorem{conjecture}{Conjecture}

\newtheorem{example}{Example}
\newtheorem{remark}{Remark}

\numberwithin{equation}{section}

\begin{document}

\sloppy

\begin{center}
{\bf \Large Local high-degree polynomial integrals of geodesic flows and the generalized hodograph method\footnote{This work is supported by the grant of the Russian Science Foundation No. 24-11-00281, https://rscf.ru/project/24-11-00281/}}
\end{center}

\medskip

\begin{center}
{\bf Sergei Agapov}
\end{center}

\medskip

\begin{quote}
\noindent{\bf Abstract.}{\small We study Riemannian metrics on 2-surfaces with integrable geodesic flows such that an additional first integral is high-degree polynomial in momenta. This problem reduces to searching for solutions to certain quasi-linear systems of PDEs which turn out to be semi-Hamiltonian. We construct plenty of local explicit and implicit integrable examples with polynomial first integrals of degrees 3, 4, 5. Our construction is essentially based on applying the generalized hodograph method.
}

\medskip


\end{quote}

\medskip

\section{Introduction and the main results}

Consider a 2-dimensional surface $M$ with local coordinates $x=(x^1,x^2)$ and a Riemannian metric $ds^2=g_{ij}(x)dx^idx^j.$ The Hamiltonian system
\begin{equation}
\label{Hamsyst}
\dot{x}^j = \{x^j,H\}, \qquad \dot{p}_j = \{p_j,H\}, \qquad H = \frac{1}{2} g^{ij} p_ip_j, \qquad i,j=1,2,
\end{equation}
with the standard Poisson bracket
$$
\{f,H\} = \sum_{i=1}^2 \left ( \frac{\partial f}{\partial x^i} \frac{\partial H}{\partial p_i} - \frac{\partial f}{\partial p_i} \frac{\partial H}{\partial x^i} \right )
$$
defines {\it the geodesic flow} of the given metric. Namely, projections of trajectories of system~\eqref{Hamsyst} from the cotangent bundle $T^*M$ to the configurational space $M$ coincide with geodesic curves.

A function $F(x, p)$ is called {\it a first integral} of the geodesic flow~\eqref{Hamsyst} if it is preserved along its trajectories, i.e. if $\frac{dF}{dt} = \{F,H\} \equiv 0.$ The geodesic flow~\eqref{Hamsyst} is called {\it completely integrable} if in addition to $H$ there exists one more first integral $F$ such that $F,$ $H$ are functionally independent almost everywhere. In this case equations~\eqref{Hamsyst} can be integrated in quadratures (\cite{Arn}).

The problem of searching for Riemannian metrics on 2-surfaces with an integrable geodesic flow is classical: it has been studied intensively during long period of time. Let us mention a famous result of Jacobi who integrated the geodesic equations on the three-axis ellipsoid in terms of elliptic functions (\cite{Jac}). Another classical integrable example is the geodesic flow on surfaces of revolution, an additional first integral (Clairaut's integral) is linear in momenta. We briefly expose other known results.

We start from the local aspect of this problem. Introduce locally the conformal coordinates
\begin{equation}
\label{conformal}
ds^2=\Lambda(x,y)(dx^2+dy^2)
\end{equation}
on $M.$ As known, it is always possible to do this in a small neighborhood of a generic point.

In general, in most known integrable cases, the first integrals have the form of polynomials in momenta. Polynomial integrals of low degrees are well-studied at large. In local conformal coordinates~\eqref{conformal} integrable examples with an additional linear or quadratic integral are known to be of the following form (\cite{Mas}---\cite{Birk1}, see also~\cite{BMF}):
\begin{equation}
\label{linear}
ds^2=f(x)(dx^2+dy^2), \qquad F_1 = p_2,
\end{equation}
\begin{equation}
\label{quadratic}
ds^2 = (f(x)+g(y))(dx^2+dy^2), \qquad F_2 = \frac{gp_1^2-fp_2^2}{f+g},
\end{equation}
here $f(x),$ $g(y)$ are arbitrary functions. In the first case~\eqref{linear} the metric does not depend on $y$ (such coordinate is called {\it cyclic}) and, consequently, the conjugate momenta $p_2$ is preserved. In the second case~\eqref{quadratic} we have {\it the Liouville metric} and the quadratic first integral.

It is known that in a typical situation the geodesic flow of a generic Riemannian metric is not polynomially integrable (\cite{KrMat}). However, the following theorem holds true.

\begin{theorem}(\cite{Kz1},\cite{Ten})
\label{T_KzT}
For any $n \geq 1, n \in \mathbb{N}$ there exists an analytic function $\Lambda(x,y)$ such that the geodesic flow of the metric~\eqref{conformal} admits an irreducible polynomial integral of the given degree $n$ with analytic (in a small neighborhood of a point $x=y=0$) coefficients.
\end{theorem}

We recall that a polynomial first integral $F$ of degree $n$ functionally independent on the Hamiltonian is called {\it irreducible} if there are no polynomial first integrals of degrees $m < n.$

Theorem~\ref{T_KzT} states that there exist plenty of local integrable examples of geodesic flows with additional non-trivial polynomial first integrals of high degrees. However, the problem of constructing such metrics and first integrals turned out to be very complicated. In general case this problem reduces to complicated quasi-linear systems of PDEs. As proved in~\cite{BM1}---\cite{BM3} (see also~\cite{Der}) these systems typically belong to the class of diagonalizable {\it hydrodynamic type systems} (\cite{DN1},~\cite{DN2}) and have many beautiful properties, for instance, they turned out to be {\it semi-Hamiltonian} (\cite{Ts1}). This fact is crucial since it allows one to apply {\it the generalized hodograph method} (\cite{Ts1}) to these systems to construct its solutions. Various attempts to do this have been made in~\cite{PT1} (see also~\cite{MP1}) where an explicit algorithm was proposed for solving this problem. Unfortunately, this algorithm turned out to be very non-trivial and produced no new integrable examples. One of significant difficulties is that for $n>2$ it is not easy to reduce an initial non-diagonal semi-Hamiltonian system to a diagonal form explicitly. As far as we know the first example of the polynomially integrable geodesic flow obtained via the generalized hodograph method was constructed in~\cite{AbdMir}, the additional polynomial first integral in this example is of the fourth degree.

In this paper we apply the generalized hodograph method combined with some other ideas. The key object in the generalized hodograph method are commuting flows (symmetries). Constructing such flows for non-diagonal semi-Hamiltonian systems is a complicated problem. We split the process of searching for these flows into two stages, one of which is purely algebraic in nature and can therefore be easily done. Our approach significantly simplifies the problem of constructing symmetries and allows us to overcome main difficulties which arose in~\cite{PT1}. Our main result is the following: we construct plenty of explicit and implicit examples of 2-dimensional Riemannian metrics with integrable geodesic flows; additional first integrals are polynomials in momenta of degrees 3, 4, 5 (see Theorems~\ref{n3maintheorem} ---~\ref{n5maintheorem} and Examples~\ref{ex1} ---~\ref{ex9} below). It is important to notice that typically examples obtained via the generalized hodograph method turn out to be implicit. However, sometimes it is possible to give them an explicit form via an appropriate canonical change of coordinates.

Another interesting problem is to investigate so-called {\it superintegrable} geodesic flows (\cite{Koenigs},~\cite{TWH}). Recall that the geodesic flow~\eqref{Hamsyst} is called {\it superintegrable} if in addition to the Hamiltonian $H$ it admits two more first integrals $F_1,$ $F_2$ (i.e. $\{F_1, H\} = \{F_2, H\} = 0$) such that all of them are functionally independent almost everywhere. In most known examples of 2-dimensional superintegrable geodesic flows one of the first integrals is linear or quadratic in momenta and another one is polynomial of a higher degree (\cite{MatvSh}---\cite{Val2}) or non-polynomial (\cite{AS1},~\cite{GS}). We note that an existence of a linear or a quadratic first integral is crucial for results obtained in~\cite{MatvSh}---\cite{Val2} . In contrast with all these results, our main aim is to construct irreducible high-degree polynomial integrals.

The global aspect of this problem is also of an undoubted interest. A special attention is usually paid to the case of compact configurational spaces. In this case certain topological obstacles to global integrability are known to exist. Namely, if the genus of a surface is larger than 1, then there exist no analytical Riemannian metrics with an additional analytical first integral independent on the Hamiltonian (\cite{Kz2}).

Integrable geodesic flows on the 2-sphere with linear and quadratic integrals were described in~\cite{Kol1}, they have rather complicated forms. Integrable examples with polynomial integrals of degree 3 and 4 are also known to exist on the 2-sphere: they can be obtained from famous integrable cases (Euler, Lagrange, Kovalevskaya cases etc.) in the dynamics of a rigid body by applying the Maupertuis principle (\cite{BKF}).

As also shown in~\cite{Kol1}, Riemannian metrics on the 2-torus with linear and quadratic integrals have exactly the forms~\eqref{linear},~\eqref{quadratic}. The following conjecture is known about polynomial integrals of higher degrees.

\begin{conjecture}(\cite{BKF})
\label{Conj}
The largest possible degree of any irreducible polynomial first integral of the geodesic flow on the 2-torus is not larger than 2.
\end{conjecture}

Although this conjecture has been proven in many special cases via various methods and approaches (\cite{KD1},~\cite{KD2},~\cite{BM2},~\cite{BM3},~\cite{T1},~\cite{Sharaf}), in general case it still remains an open problem.

Given a high-degree polynomial integral of the geodesic flow, it is an important and usually very non-trivial problem to verify the absence of polynomial integrals of lesser degrees. One of possible approaches to this problem based on the theory of differential invariants was developed in~\cite{Krugl}. An explicit criterion for existence of linear integrals described in~\cite{Krugl} is the following.

\begin{theorem}(\cite{Krugl})
\label{Criteria}
Given a 2-dimensional metric $ds^2=g_{ij}du^idu^j,$ consider the following functions:

- $R = R^i_{ijk}g^{jk}$ (scalar curvature), where $R^i_{ljk}$ is the curvature tensor of the given metric;

- $L = g^{ij} \frac{\partial R}{\partial u^i} \frac{\partial R}{\partial u^j};$

- $\Delta = \frac{1}{\sqrt{det \ g}} \frac{\partial}{\partial u^i} \left( g^{ij} \sqrt{det \ g} \frac{\partial R}{\partial u^j} \right).$

If the geodesic flow of the given metric admits a linear in momenta first integral, then $R$ and $L$ as well as $R$ and $\Delta$ are functionally dependent, i.e.
\begin{equation}
\label{Criteria_formulae}
det \begin{pmatrix}
	\frac{\partial R}{\partial u^1} & \frac{\partial R}{\partial u^2}  \\
	\frac{\partial L}{\partial u^1} & \frac{\partial L}{\partial u^2}
\end{pmatrix} = det \begin{pmatrix}
	\frac{\partial R}{\partial u^1} & \frac{\partial R}{\partial u^2}  \\
	\frac{\partial \Delta}{\partial u^1} & \frac{\partial \Delta}{\partial u^2}
\end{pmatrix} = 0.
\end{equation}

\end{theorem}
In one form or another this criterion has also been described and/or used in many other papers (e.g., see~\cite{Darboux1}---\cite{Matv1}).
We shall use this criterion to verify the absence of linear integrals in all of our constructed examples below.

A natural generalization of polynomial integrals are the ones of the form
\begin{equation}
\label{generalintegral}
F = \prod_{j=1}^N (a_j(x,y)p_1+b_j(x,y)p_2)^{n_j}, \qquad n_j \in \mathbb{R}.
\end{equation}
Such integrals were firstly studied by Darboux in~\cite{Darboux1} and later in many papers (e.g., see~\cite{MP1},~\cite{AS2},~\cite{DemI}). Depending on $n_j,$ these integrals may be polynomial, rational, algebraic or even transcendental functions. In the case $N=2,$ $n_1=1$ explicit examples of such integrals (with arbitrary $n_2=n \in \mathbb{R}$) were constructed in~\cite{AS2}. We provide here one of these examples which is expressed in terms of elementary functions. We note that there are some misprints in~\cite{AS2} and the correct example is the following.

\begin{example} (\cite{AS2})
\label{ex0}
Suppose that $n(n+1)(n+2) \neq 0.$ Then the geodesic flow of the metric
$$
ds^2 = e^{\frac{4(n+1)}{n}x} \{ (2n^2 \sin^2((n+1)y)+8\cos^2((n+1)y))dx^2
$$
$$
+ 2n(n+1)(n-2)\sin(2(n+1)y)dxdy+ n^2(n^2-2n+2+n(n-2)\cos(2(n+1)y))  dy^2 \}
$$
admits an additional first integral of the form
$$
F = e^{-\frac{n^2+3n+2}{n}x} \left(3n-2+(n-2)\cos(2(n+1)y) \right)^{-1-n} \left( np_1\sin((n+1)y)+2p_2\cos((n+1)y) \right)^n  \times
$$
$$
\times \left( (n-2) (np_1\sin(2(n+1)y)+p_2\cos(2(n+1)y)) - (n+2)p_2 \right).
$$
If in addition $n \neq 2,$ then this metric is not flat.

\end{example}

We note that for any natural $n \in \mathbb{N},$ $n>2$ the additional integral in Example~\ref{ex0} is polynomial in momenta of degree $n+1$ and there are no linear integrals (it can be verified via Theorem~\ref{Criteria}).

The simplest non-polynomial integrals~\eqref{generalintegral} are the ones having a rational form. We refer the reader to~\cite{Kz3},~\cite{AS1},~\cite{AS2},~\cite{AD1} where many examples of rational integrals were constructed.

This paper is organized as follows. In section 2 we give the definitions of semi-Hamiltonian systems, their commuting flows (symmetries) and briefly describe the generalized hodograph method (\cite{Ts1}). In section 3 we recall the results obtained in~\cite{BM1} about semi-Hamiltonicity of quasi-linear systems of PDEs arising in the problem of polynomially integrable geodesic flows on the 2-torus. In section 4 we describe well-known cases of integrals of small degrees $n=1, 2.$ In section 5 we explain how to search for commuting flows of non-diagonal semi-Hamiltonian systems. Following this approach, in sections 6 --- 8 we construct plenty of implicit and explicit integrable examples of metrics with polynomial integrals of degrees 3, 4 and 5.

\section{Semi-Hamiltonian systems and the generalized hodograph method}

In this section we briefly describe some fundamental properties of semi-Hamiltonian systems and recall the notions of the classical and the generalized hodograph methods.

First consider a homogeneous quasi-linear system of two differential equations on two unknown functions $u(t,x),$ $v(t,x):$
\begin{equation}
\label{classical_hod1}
u_t+a_{11}u_x+a_{12}v_x=0, \qquad v_t+a_{21}u_x+a_{22}v_x=0,
\end{equation}
where $a_{ij} = a_{ij}(u,v)$ for any $i, j = 1, 2.$ Let us swap independent and dependent variables, i.e. we assume that $x=x(u,v),$ $t=t(u,v).$ It is easy to verify that if $u_xv_t-u_tv_x \neq 0,$ then quasi-linear system~\eqref{classical_hod1} is equivalent to the linear one:
\begin{equation}
\label{classical_hod2}
x_v=a_{11}t_v-a_{12}t_u, \qquad -x_u=a_{21}t_v-a_{22}t_u.
\end{equation}
This is how {\it the classical hodograph method} works (e.g., see~\cite{RYa}).

\vspace{2mm}

The generalized hodograph method (\cite{Ts1}) is a natural generalization of the classical one.

A diagonal quasi-linear system of equations (no summation over $i$)
\begin{equation}
\label{semiHamsyst}
r^{i}_t = v_i(r)r^{i}_x, \qquad i = 1,...,n, \qquad v_i \neq v_j
\end{equation}
is called {\it semi-Hamiltonian} if its coefficients $v_i(r)$ are distinct in a domain and for $n>2$ satisfy
\begin{equation}
\label{semiHamcond}
\partial_i\left(\frac{\partial_j v_k}{v_j - v_k}\right) = \partial_j\left(\frac{\partial_i v_k}{v_i - v_k}\right),\quad i \neq j \neq k \neq i.
\end{equation}
In the case $n\leq2$ any hyperbolic system is semi-Hamiltonian.

The evolutionary system
\begin{equation}
\label{symm}
r^{i}_\tau = w_i(r)r^{i}_x, \qquad i = 1,...,n
\end{equation}
is called {\it a symmetry of hydrodynamic type} (or, equivalently, {\it a commuting flow}) of~\eqref{semiHamsyst} if and only if $\partial_\tau(r^{i}_t) = \partial_t(r^{i}_{\tau})$ for any $i = 1,...,n.$ Components $w_i(r)$ of the symmetries satisfy the following relations:
\begin{equation}
\label{symm_cond}
\frac{\partial_k v_i}{v_k-v_i} = \frac{\partial_k w_i}{w_k-w_i} , \qquad i \neq k.
\end{equation}

\begin{theorem}(\cite{Ts1})
\label{T_symm}
A semi-Hamiltonian diagonal system~\eqref{semiHamsyst} has infinitely many symmetries of hydrodynamic type: the flows~\eqref{symm} commuting with it, that can be parameterized by $n$ functions of one variable. All these flows commute with each other and their matrices are diagonal.
\end{theorem}

Consider the following system of $n$ equations on $n$ unknown functions $r^i:$
\begin{equation}
\label{gen_hod}
w_i(r) = v_i(r) t + x, \quad i = 1,...,n.
\end{equation}

\begin{theorem}(\cite{Ts1})
\label{T_hod}
A smooth solution $r^i(t,x)$ of system~\eqref{gen_hod} is a solution of the diagonal semi-Hamiltonian system~\eqref{semiHamsyst}. Conversely, any solution $r^i(t,x)$ of~\eqref{semiHamsyst} can be locally represented as a solution of~\eqref{gen_hod} in a neighborhood of a point $(t_0, x_0)$ such that $r^i_x(t_0, x_0)\neq0$ for every $i,$ where $w_i(u)$ are coefficients of a hydrodynamic flow~\eqref{symm} commuting with~\eqref{semiHamsyst}.
\end{theorem}

In the case $n=2$ the generalized hodograph method is equivalent to the classical one.

\vspace{2mm}

We note that typically semi-Hamiltonian systems arising in applications are written in a non-diagonal form:
\begin{equation}
\label{semiHamnondiag}
u^{i}_t = \sum_{j=1}^{n} v^{i}_j(u)u^{j}_x,\quad i = 1,..., n.
\end{equation}
In this case the generalized hodograph method can be applied in the following way.
Let us search for commuting flows in the form
\begin{equation}
\label{nondiagsymm}
u^{i}_\tau = \sum_{j=1}^{n} w^{i}_j(u)u^{j}_x,\quad i = 1,..., n
\end{equation}
such that the following conditions hold true:
\begin{equation}
\label{nondiagcompcond}
\partial_\tau(u^{i}_t) = \partial_\tau \left(\sum_{j = 1}^{n} v^{i}_j(u)u^{j}_x\right) = \partial_t \left(\sum_{j = 1}^{n} w^{i}_j(u)u^{j}_x\right) = \partial_t(u^{i}_\tau).
\end{equation}
Consider then the following system
\begin{equation}
\label{nondiaghodog}
x \delta^{i}_k + t v^{i}_k = w^{i}_k.
\end{equation}
This system contains $n^2$ equations on $n$ unknown functions $u^i.$ However, as proved in~\cite{Ts1}, if system~\eqref{semiHamnondiag} is semi-Hamiltonian and $w^{i}_k$ define its symmetry (i.e. if relations~\eqref{nondiagcompcond} hold true), then system~\eqref{nondiaghodog} must be compatible. Having solved~\eqref{nondiaghodog} w.r.t. $u^i(t,x),$ we automatically obtain solutions to system~\eqref{semiHamnondiag}.

\section{Polynomial integrals of geodesic flows on the 2-torus and semi-Hamiltonian systems}

In this section we describe a connection between polynomial integrals of geodesic flows and semi-Hamiltonian systems which was discovered in~\cite{BM1}.

\begin{theorem}(\cite{BM1})
\label{T_BM1}
Suppose that Hamiltonian system~\eqref{Hamsyst} on the two-torus admits an additional integral $F$ which is a homogeneous polynomial of degree $n.$ Then on the covering plane $\mathbb{R}^2$ there exist the global coordinates $(t,x)$ such that the metric has the following form:
$$ds^2=g^2(t,x)dt^2+dx^2,$$
and the integral $F$ can be written in the form
$$F=\displaystyle\sum_{k=0}^{n} \frac{a_k(t,x)}{g^{n-k}}p_1^{n-k}p_2^{k},$$
where the last two coefficients can be normalized in the following way: $a_{n-1}\equiv g,$ $a_{n} \equiv 1.$ Then the commutation relation $\left\{F,H\right\}=0$ is equivalent to the system of $n$ quasi-linear equations on the unknown functions $a_{0},...,a_{n-1}$ of the form
\begin{equation}
\label{key_system_general}
u_t^i+v_j^i(u)u_x^j=0,
\end{equation}
where $u^i=(a_{0},...,a_{n-1})^T$ and the matrix $v_j^i$ has the form:
\begin{equation}
\label{key_matrix_general}
v_j^i =\begin{pmatrix}
		0 & 0 & \cdots & 0 & 0 & a_1 \\
		a_{n-1} & 0 & \cdots & 0 & 0 & 2a_{2}-na_{0} \\
		\cdots  & \cdots  & \cdots & \cdots & \cdots & \cdots  \\
		0 & 0 & \cdots & a_{n-1} & 0 & \big(n-1\big)a_{n-1}-3a_{n-3}\\
		0 & 0 & \cdots & 0 & a_{n-1} & na_{n}-2a_{n-2}
\end{pmatrix}.
\end{equation}
Functions $a_i$, $g$ are periodic w.r.t. $x$ and quasi-periodic w.r.t. $t.$
\end{theorem}

\begin{theorem}(\cite{BM1})
\label{T_BM2}
System~\eqref{key_system_general} has the following properties.

1. In the region of hyperbolicity (all eigenvalues are real and distinct) there
exists a change of variables (Riemann invariants) $(a_0,\dotsc, a_{n-1})\to\mathbb(r_1,\dotsc, r_n)$ which transforms system~\eqref{key_system_general} to a diagonal form:
$$(r_i)_t+\lambda_i(r_1, \dotsc ,r_n)(r_i)_x=0, \qquad i=1, \dotsc ,n;$$

2. there exists a regular change of variables $(a_0,\dotsc, a_{n-1})\to\mathbb(G_1,\dotsc, G_n)$ such that system~\eqref{key_system_general} can be written in the form of conservation laws: $$(G_i(a_0, \dotsc, a_{n-1}))_t+(H_i(a_0, \dotsc, a_{n-1}))_x=0, \qquad i=1, \dotsc ,n.$$
\end{theorem}

It turns out that quasi-linear systems having such a pair of properties as in Theorem~\ref{T_BM2} are necessarily semi-Hamiltonian. More precisely, the following theorem holds true.

\begin{theorem}(\cite{Sevennec},~\cite{Serre})
\label{T_symm2}
The hyperbolic diagonal system~\eqref{semiHamsyst} can be written in the form of conservation laws $$
g_i(r)_t+h_i(r)_x=0, \qquad i=1, \dotsc ,n
$$
if and only if condition~\eqref{semiHamcond} holds true.
\end{theorem}

This statement is crucial since it guaranties that system~\eqref{key_system_general} arising in the problem of integrable geodesic flows is in fact semi-Hamiltonian. This allows to apply the generalized hodograph method to construct solutions to~\eqref{key_system_general}.

Periodic solutions to system~\eqref{key_system_general} are known to exist in cases $n=1, 2$ (see section 4 below). The question about existence of nontrivial periodic solutions in the case $n > 2$ remains open (see Conjecture~\ref{Conj} in Introduction).

On the other hand, in the case $n > 2$ system~\eqref{key_system_general} admits nontrivial non-periodic solutions which give rise to local high-degree polynomial integrals of geodesic flows. The main aim of this paper is to describe these solutions. We notice that though our results do not solve Conjecture~\ref{Conj}, they still may be useful for further investigations of this problem.

\section{Cases $n=1,$ $n=2$}

In the simplest cases $n=1$ and $n=2$ solutions to system~\eqref{key_system_general} are well-known (e.g., see~\cite{PT1}). In this section we expose them for the sake of completeness.

The case $n=1$ is trivial: the general solution is given by
$$a_0(t,x) = g(t,x) = f(t-x), \qquad a_1(t,x)\equiv 1,$$
here $f$ is an arbitrary function of one argument. So in this case we obtain
$$
H=\frac{1}{2}\left(\frac{p_1^{2}}{f^2(t-x)}+p_2^{2}\right), \qquad F=p_1+p_2, \qquad \left\{F,H\right\}=0.
$$

\vspace{2mm}

Now let us demonstrate how the generalized hodograph method can be applied to system~\eqref{key_system_general} in the case $n=2.$ A general solution in this case can be obtained in an implicit form.
We have
$$ds^2=g^2(t,x)dt^2+dx^2, \quad H=\frac{1}{2}\left(\frac{p_1^{2}}{g^2(t,x)}+p_2^{2}\right), \quad F=\frac{a_0(t,x)}{g^2}p_1^2+\frac{a_1(t,x)}{g}p_1 p_2+a_2(t,x)p_2^2.
$$
Due to $a_1(t,x)\equiv g(t,x),$ $a_2(t,x)\equiv 1$ we obtain that the condition $\left\{F,H\right\}=0$ implies
$$
(a_0)_t+gg_x=0, \qquad g_t+2(1-a_0)g_x+g(a_0)_x=0.
$$
This system can be written in the form of conservation laws:
$$
\left(a_0\right)_t+\left( \frac{g^2}{2} \right)_x=0, \qquad	\left( \dfrac{1}{2g^2} \right)_t+ \left( \frac{1-a_0}{g^2} \right)_x=0.
$$
In Riemann invariants $r^1,$ $r^2:$
$$
a_0(t,x)=1-r^1(t,x)-r^2(t,x), \qquad g^2(t,x)=-4 r^1(t,x) r^2(t,x),
$$
the system takes the diagonal form:
\begin{equation}
\label{n2diag}
\left(
\begin{array}{c}
	r^1 \\
	r^2
\end{array}
\right)_t+\begin{pmatrix}
	2r^2 & 0  \\
	0 & 2r^1
\end{pmatrix}\left(
\begin{array}{c}
	r^1 \\
	r^2
\end{array}\right)_x=0,
\end{equation}
i.e. $r_t^i+v_i(r) r_x^i=0,$ where $v_1=2r^2,$ $v_2=2r^1$. Notice that system~\eqref{n2diag} has the following property:
\begin{equation}
\label{n2weaknonl}
\frac{\partial v_1}{\partial r^1} = \frac{\partial v_2}{\partial r^2} = 0.
\end{equation}
Diagonal systems satisfying condition~\eqref{n2weaknonl} are called {\it weakly nonlinear} (\cite{RYa}). Diagonal weakly nonlinear semi-Hamiltonian systems were described in~\cite{Fer1} (see also~\cite{Pavlov},~\cite{Fer2}), various methods for constructing its solutions were also discussed there.

Weakly nonlinear systems have the following important property: as shown in~\cite{RS1}, the absence of a gradient catastrophe is typical for solutions to such systems. In more details, for an arbitrary diagonal weakly nonlinear system
\begin{equation}
\label{weaknonlsyst}
r_t^i+\xi_i(r) r_x^i=0, \qquad \frac{\partial \xi_i}{\partial r^i}=0, \qquad i = 1, \ldots, n
\end{equation}
consider the Cauchy problem $r^i(0,x) = r^i_0(x).$ Assume that a solution exists and is bounded for $0 \leq t \leq T,$ the functions $\xi_i$ and their first derivatives are bounded in the solution and system~\eqref{weaknonlsyst} is hyperbolic in the narrow sense, i.e.
\begin{equation*}
\label{narrowhyperbolic}
\left| \xi_k(t,x,r(t,x)) \right| \leq M, \ \left| \xi_i(t,x,r(t,x)) - \xi_j(t,x,r(t,x)) \right| \geq \varepsilon > 0, \ i,j,k=1, \ldots, n, \ i \neq j,
\end{equation*}
here $M,$ $\varepsilon$ are constants. In addition assume that at the moment $t=0$ the derivatives are bounded:
\begin{equation*}
\label{boundderivatives}
\left| \frac{dr^i_0(x)}{dx} \right| \leq Q,
\end{equation*}
where $Q$ is a constant. Then, as proved in~\cite{RS1}, derivatives $\partial r^i/\partial x$ of a solution to~\eqref{weaknonlsyst} are bounded for finite $t.$

This implies that system~\eqref{n2diag} and, consequently, system~\eqref{key_system_general} in the case $n=2$ admits smooth solutions. It is interesting that as shown in~\cite{AF1}, in the case $n>2$ system~\eqref{key_system_general} is no more weakly nonlinear and this is one of the reasons why Conjecture~\ref{Conj} seems to be true.

Now we apply the generalized hodograph method to construct solutions to system~\eqref{n2diag}. Commuting flows have the form
$$
\left(
\begin{array}{c}
	r^1 \\
	r^2
\end{array}
\right)_{\tau}=\begin{pmatrix}
	w_1 & 0  \\
	0 & w_2
\end{pmatrix}\left(
\begin{array}{c}
	r^1 \\
	r^2
\end{array}\right)_x,
$$
where $w_1(r),$ $w_2(r)$ are unknown functions so far. Relations~\eqref{symm_cond} take the form
$$
\frac{\partial w_1}{\partial r^2} = \frac{w_1-w_2}{r^2-r^1}, \qquad \frac{\partial w_2}{\partial r^1} = \frac{w_2-w_1}{r^1-r^2},
$$
which implies $\partial w_1/\partial r^2=\partial w_2/\partial r^1.$
Introduce a new function $\Psi(r^1, r^2)$ such that $\Psi_{r^1}=w_1,$ $\Psi_{r^2}=w_2.$
Then $\Psi$ satisfies the Euler-Poisson-Darboux equation:
\begin{equation}
\label{EPDeq}
\Psi_{r^1 r^2}+\dfrac{\Psi_{r^1}-\Psi_{r^2}}{{r^1-r^2}}=0.
\end{equation}
The general solution to this equation has the form (e.g., see~\cite{Tricomi}):
$$
\Psi(r^1, r^2)=2u(r^1)+2v(r^2)+(r^1-r^2)(v'(r^2)-u'(r^1)),
$$
where $u(r^1),$ $v(r^2)$ are two arbitrary functions of one argument. Finding $w_1,$ $w_2$ and substituting them into~\eqref{gen_hod} we obtain the general solution to system~\eqref{n2diag} (and also to system~\eqref{key_system_general} for $n=2$) in the implicit form:
$$
t=-\frac{1}{2}(u''(r^1)+v''(r^2)), \qquad x=u'(r^1)+v'(r^2)-r^1u''(r^1)-r^2v''(r^2).
$$

Notice that the obtained solution with the quadratic first integral corresponds exactly to the Liouville metric in conformal coordinates: the corresponding change of coordinates is given in~\cite{PT1}.

\section{Description of commuting flows of non-diagonal semi-Hamiltonian systems}

In the case $n = 2,$ for the system written in Riemann invariants~\eqref{n2diag} commuting flows were completely described in section 4 which allowed to construct the general solution via the generalized hodograph method. For $n>2,$ the search for Riemann invariants of a semi-Hamiltonian system in general becomes a very complicated problem. Therefore, in this case, for constructing solutions by the generalized hodograph method, it is reasonable to describe commuting flows~\eqref{nondiagsymm} of the initial non-diagonal system~\eqref{semiHamnondiag} without its explicit diagonalization.

Consider non-diagonal semi-Hamiltonian system~\eqref{semiHamnondiag} and its commuting flow~\eqref{nondiagsymm}. Denote $U=(u^1, \ldots, u^n)^T,$ $V = (v^i_j),$ $W=(w^i_j).$

\begin{lemma}(\cite{Ts1})
\label{nondiagsymmcondition}
Two flows
\begin{equation}
\label{matrixformsymm}
U_t=V(U)U_x, \qquad U_{\tau}=W(U)U_x,
\end{equation}
commute if and only if
\begin{equation}
\label{matrixformsymmcond1}
[V,W]=VW-WV \equiv 0,
\end{equation}
\begin{equation}
\label{matrixformsymmcond2}
\left( V_{\tau}-W_t+VW_x-WV_x \right) U_x \equiv 0.
\end{equation}
\end{lemma}

The proof follows immediately from the following calculation:
$$
0= (VU_x)_{\tau}-(WU_x)_t=V_{\tau}U_x+V(WU_x)_x-W_tU_x-W(VU_x)_x
$$
$$
=(VW-WV)U_{xx}+(V_{\tau}-W_t+VW_x-WV_x)U_x.
$$
To sum it up, in order to describe symmetries of a given non-diagonal semi-Hamiltonian system $U_t=V(U)U_x,$ we have to do the following steps:

1) to find all matrices $W,$ which commute with $V$ (condition~\eqref{matrixformsymmcond1}),

2) to satisfy condition~\eqref{matrixformsymmcond2}.

First condition~\eqref{matrixformsymmcond1} is purely algebraic and in principle can be immediately solved for any given semi-Hamiltonian system. Second condition~\eqref{matrixformsymmcond2} is much more complicated and in general case causes great difficulties.

Let us demonstrate how this algorithm can be applied to system~\eqref{key_system_general} in the case $n=2.$ We have
$$V =
\begin{pmatrix}
	0 & a_1 \\
	a_1 & 2-2a_0
\end{pmatrix} \qquad W =
\begin{pmatrix}
	w_1 & w_2 \\
	w_3 & w_4
\end{pmatrix},$$
where components $w_j$ of symmetries are unknown functions which depend on $a_0,$ $a_1.$
Condition~\eqref{matrixformsymmcond1} implies
$$
w_3=w_2, \qquad w_4=w_1+\frac{2(1-a_0) w_2}{a_1}.
$$
By direct calculations we obtain that condition~\eqref{matrixformsymmcond2} is equivalent to:
$$
(w_1)_{a_1}-(w_2)_{a_0}=0,  \qquad w_2-a_1 (w_2)_{a_1}+a_1 (w_1)_{a_0}-2(-1+a_0) (w_2)_{a_0}=0.
$$
Introduce a new function $\Psi(a_0,a_1)$ such that $w_1=\Psi_{a_0},$ $w_2=\Psi_{a_1},$ then $\Psi(a_0,a_1)$ satisfies the following equation:
$$
a_1 \Psi_{a_0 a_0}-2(a_0-1)\Psi_{a_0 a_1}-a_1 \Psi_{a_1 a_1}+\Psi_{a_1}=0.
$$
This equation appeared and was solved in~\cite{AF1} by the standard method of characteristics in the following way. Namely, the change of variables
$$
r_1 = 1-a_0-\sqrt{(1-a_0)^2+a_1^2}, \qquad r_2 = 1-a_0+\sqrt{(1-a_0)^2+a_1^2}
$$ reduces this equation to the canonical form~\eqref{EPDeq}, its general solution was described above in section 4. To sum it up, all the components $w_j$ of the symmetries are found. Substituting them into~\eqref{nondiaghodog} we obtain the general solution to~\eqref{key_system_general} which is obviously equivalent to the one obtained in section 4 (see all the details in~\cite{AF1}).

Below we shall apply this algorithm for describing symmetries and constructing solutions to~\eqref{key_system_general} for $n>2.$

\section{Case $n=3$}

Now we are in position to deal with system~\eqref{key_system_general} in the cases $n>2.$ In this section we consider the case $n=3.$ We have $U=(a_0, a_1, a_2)^T,$ $a_3(t,x)\equiv 1,$ $g(t,x) = a_2(t,x),$
\begin{equation}
\label{genFn3}
F=\frac{a_0}{a_2^3}p_1^3+\frac{a_1}{a_2^2}p_1^2p_2+p_1p_2^2+p_2^3
\end{equation}
and system~\eqref{key_system_general} is of the form
\begin{equation}
\label{n3system}
U_t+VU_x=0, \qquad V=
\begin{pmatrix}
		0 & 0 & a_1 \\
		a_{2} & 0 & 2a_{2}-3a_{0} \\
		0 & a_{2} & 3-2a_1
\end{pmatrix}.
\end{equation}

In this section we apply the generalized hodograph method to system~\eqref{n3system} and construct its solutions. First of all we have to describe symmetries
$$
U_\tau+W(U)U_x=0
$$
of~\eqref{n3system}. Due to Lemma~\ref{nondiagsymmcondition} conditions~\eqref{matrixformsymmcond1},~\eqref{matrixformsymmcond2} must be fulfilled. We start with searching for matrices $W$ such that $[V,W]=0.$ This condition has the form of a linear system of equations of rank 6 which consists of 9 equations on the unknown components of $W.$ Solving it we obtain

\begin{lemma}
\label{n3lemma}
For matrix $V$ of the form~\eqref{n3system} the condition $[V,W]=0$ holds true if and only if $W$ has the form
\begin{equation}
\label{n3symm}
W=
\begin{pmatrix}
		(3a_0-2a_2)P+(2a_1-3)R+S & a_1P & a_1R \\
		(2a_1-3)P+a_2R & (2a_1-3)R+S & a_1P+(2a_2-3a_0)R \\
		a_2P & a_2R & S
\end{pmatrix}
\end{equation}
with arbitrary $P,$ $R,$ $S.$
\end{lemma}

Now the first condition~\eqref{matrixformsymmcond1} is fulfilled and we need to fulfill the second one~\eqref{matrixformsymmcond2}. The following lemma can be proved by direct calculations.
\begin{lemma}
\label{n3lemma_add1}
Condition~\eqref{matrixformsymmcond2} holds true if and only if the unknown functions $P(a_0,a_1,a_2),$ $R(a_0,a_1,a_2),$ $S(a_0,a_1,a_2)$ satisfy the following system of PDEs:
\begin{equation}
\label{n3pde}
\begin{gathered}
a_1 \partial_0 P = (3a_0-2a_2) \partial_0 R + (2a_1-3) \partial_1 R + a_2 \partial_2 R,\hfill\\
\partial_1 P = \partial_0 R, \qquad \partial_2 P = \partial_1 R,\hfill\\
\partial_0 S = a_2 \partial_1 R - 2P, \qquad \partial_1 S = a_2 \partial_2 R - 2R,\hfill\\
\partial_2 S = a_1 \partial_0 R + (2a_2-3a_0) \partial_1 R + (3-2a_1) \partial_2 R + P,\hfill
\end{gathered}
\end{equation}
where $\partial_i = \frac{\partial}{\partial {a_i}},$ $i=0,1,2.$
\end{lemma}

Lemmas~\ref{n3lemma},~\ref{n3lemma_add1} give the complete description of symmetries of system~\eqref{n3system}. According to the generalized hodograph method, solutions to~\eqref{n3system} are given by system~\eqref{nondiaghodog}.
To sum it up, the following theorem describes solutions to~\eqref{n3system}.

\begin{theorem}
\label{n3maintheorem}
Assume that functions $P(a_0,a_1,a_2),$ $R(a_0,a_1,a_2),$ $S(a_0,a_1,a_2)$ satisfy system~\eqref{n3pde}. Then any smooth solution $a_k(t,x)$ to the system of equations
\begin{equation}
\label{n3eq}
P=0, \qquad R=t, \qquad S=(3-2a_1)t-x
\end{equation}
is also a solution to system~\eqref{n3system}.
\end{theorem}

Let us prove Theorem~\ref{n3maintheorem}. Lemmas~\ref{n3lemma},~\ref{n3lemma_add1} imply that matrix $W$ of the form~\eqref{n3symm} with functions $P,$ $R,$ $S$ satisfying~\eqref{n3pde} defines the symmetry of~\eqref{n3system}. It is left to notice that system~\eqref{nondiaghodog} with such $W$ is equivalent to relations~\eqref{n3eq}. This completes the proof.

Let us make one remark concerning the statement of Theorem~\ref{n3maintheorem}.

\begin{remark}
\label{rem1}

In order to construct solutions to~\eqref{n3system} via Theorem~\ref{n3maintheorem}, we have to be able to solve~\eqref{n3eq} w.r.t. $a_0,$ $a_1,$ $a_2$ as functions depending on $t,$ $x.$ For instance, functions $P(a_0,a_1,a_2) = R(a_0,a_1,a_2) = S(a_0,a_1,a_2) \equiv 0$ obviously satisfy~\eqref{n3pde} but they do not produce any solutions to~\eqref{n3system} since in this case equations~\eqref{n3eq} cannot be solved w.r.t. $a_k(t,x),$ $k=0,1,2.$

\end{remark}

\vspace{3mm}

It is not easy to construct a general solution to~\eqref{n3pde}. However, plenty of partial solutions to this system can be obtained. The simplest way to do this seems to be the following. Suppose that one of the unknown functions has a certain simple form. For instance, assume that $P$ is a polynomial in $a_k$ of the first degree with unknown constant coefficients, that is $$P(a_0, a_1, a_2) = k_0a_0+k_1a_1+k_2a_2+k_3, \qquad k_j \in \mathbb{R}.$$ Let us substitute such $P(a_0, a_1, a_2)$ into system~\eqref{n3pde}. Integrating this system we find the unknown functions $$P=(2a_0+a_2)k_2+k_3, \qquad R=(a_1+3\log a_2)k_2+k_4,$$ $$S=3a_2^2k_2/2-a_0(2k_3+(2a_0+a_2)k_2)+(3-a_1)a_1k_2-2a_1k_4+a_2k_3+(9-6a_1)k_2\log a_2+k_5,$$ here $k_2, \ldots, k_5$ are arbitrary constants. System~\eqref{n3eq} with such functions $P,$ $R,$ $S$ defines a family of implicit solutions to~\eqref{n3system} parameterized by several arbitrary constants. Assuming for simplicity $k_3=k_4=k_5=0$ and denoting $k_2=k$ we obtain the following

\begin{example}
\label{ex1}
Assume that functions $a_0(t,x),$ $a_1(t,x),$ $a_2(t,x)$ satisfy
\begin{equation}
\label{ex1impl1}
2a_0+a_2=0, \qquad k(a_1+3 \log a_2)=t, \qquad 2x+k(2a_1^2+3a_2^2)=0,
\end{equation}
here $k\in\mathbb{R}$ is an arbitrary constant.

Then the geodesic flow of the metric $ds^2 = g^2 (t,x) dt^2 + dx^2$ (recall that $a_2\equiv g$) admits the cubic in momenta first integral $F_3$ of the form
$$
F_3 = \frac{a_0}{g^3}p_1^{3}+\frac{a_1}{g^{2}}p_1^{2}p_2+p_1p_2^2+p_2^3.
$$
\end{example}

Solving the first two equations of~\eqref{ex1impl1} we obtain that
$$
a_0=-a_2/2, \qquad a_1 = -3 \log a_2 +t/k
$$
and the metric coefficient $a_2(t,x)$ satisfies the transcendental equation:
$$
3ka_2^2-12t\log a_2+18k\log^2a_2 + 2x + 2t^2/k = 0.
$$
In order to get rid of this implicitness, let us go back to~\eqref{ex1impl1}. We have $a_0=-a_2/2$ and
\begin{equation}
\label{ex1impl2}
k(a_1+3 \log a_2)=t, \qquad 2x+k(2a_1^2+3a_2^2)=0.
\end{equation}
Let us make the change of coordinates $(t,x) \rightarrow (a_1,a_2)$ considering~\eqref{ex1impl2} as the corresponding formulae. New conjugate momenta $P_1,$ $P_2$ can be found from the relation
$$
p_1dt+p_2dx=P_1da_1+P_2da_2.
$$
Without loss of generality assume that $k=1.$ Rewriting Hamiltonian and the first integral in terms of new coordinates and momenta we finally obtain the following explicit integrable example.

\begin{example}
\label{ex2}
Denote $u^1=x,$ $u^2=y.$ The geodesic flow of the metric $ds^2 = g_{ij}(u)du^idu^j$ where
$$
g_{11} = 4x^2+y^2, \qquad g_{12} = 3y(1+2x), \qquad g_{22} = 9\left(1+y^2\right)
$$
admits the cubic in momenta first integral $F=\alpha f_{k}(x,y)p_1^{3-k}p_2^k,$ where
$$
f_0 = 27 \left( y^4-2y^2(1+x)-2 \right), \qquad f_1 = 18y \left(4x^2-2y^2(x-1)+2x+3 \right),
$$
$$
f_2 = -6 \left( y^4+4x^3-y^2 (2x^2-4x-3) \right), \qquad f_3 = 2y^3 \left(2x+1 \right),
$$
here $\alpha = \left(y^2-2x\right)^{-3}.$ Gauss curvature $K$ of this metric is equal to
$$
K = \frac{y^2+2x+6}{9\left(y^2-2x\right)^3}.
$$
\end{example}

We note that the geodesic flow of the metric constructed in Example~\ref{ex2} does not admit any linear in momenta first integrals. This can be proved with the help of Theorem~\ref{Criteria}. This also implies the absence of linear integrals in Example~\ref{ex1}. It means that the constructed cubic integrals in these examples cannot be expressed in terms of polynomial integrals of lesser degrees.

Now assume that $P$ is a non-homogeneous quadratic polynomial in $a_0,$ $a_1,$ $a_2$ with unknown constant coefficients. As previously, by substituting such $P$ into system~\eqref{n3pde} and integrating it, we find all the unknown functions $P(a_0, a_1, a_2),$ $R(a_0, a_1, a_2),$ $S(a_0, a_1, a_2).$ Substituting these functions into system~\eqref{n3eq} we obtain

\begin{example}
\label{ex3}
Assume that functions $a_0(t,x),$ $a_1(t,x),$ $a_2(t,x)$ satisfy
\begin{equation}
\label{ex3impl1}
\begin{gathered}
k_1(5a_0^2+a_1^2+a_2^2+2a_0a_2)+k_2(2a_0+a_2)+2k_1a_1+k_3=0,\hfill\\
2k_1a_1(a_0+a_2)+2k_1(a_0+5a_2)+k_2a_1+k_4+3k_2 \log a_2=t,\hfill\\
k_1(-10a_0^3+5a_2^3-3a_0^2a_2+6a_0a_1^2+9a_1^2a_2)-6k_2a_0^2+3k_2a_1^2\hfill\\
+9k_2a_2^2/2-18k_1a_0a_1-3k_2a_0a_2+12k_1a_1a_2\hfill\\
-6(3k_1+k_3)a_0+3k_3a_2-9k_4+3k_5=-3x,\hfill
\end{gathered}
\end{equation}
here $k_j\in\mathbb{R}$ are arbitrary constants.

Then the geodesic flow of the metric $ds^2 = g^2 (t,x) dt^2 + dx^2$ (recall that $a_2\equiv g$) admits the cubic in momenta first integral $F$ of the form
$$
F_3 = \frac{a_0}{g^3}p_1^{3}+\frac{a_1}{g^{2}}p_1^{2}p_2+p_1p_2^2+p_2^3.
$$
\end{example}

To sum it up, Example~\ref{ex1} corresponds to the linear function $P(a_0, a_1, a_2),$ Example~\ref{ex3} --- to the quadratic one. Many other integrable examples can be constructed by choosing $P$ in a polynomial form of higher degrees. They have more cumbersome forms so we do not provide them here.

Another class of integrable examples can be obtained in the following way. It follows from~\eqref{n3pde} that function $P$ satisfies the following PDE of the second order:
$$
\partial_{02} P = \partial_{11} P.
$$
Applying the method of separation of variables to this equation, namely, searching $P$ in the form $P(a_0,a_1,a_2) = \widetilde{P}(a_0) \widehat{P}(a_1,a_2)$ and integrating system~\eqref{n3pde} we find all the unknown functions $P,$ $R,$ $S.$ This allows to obtain another integrable example.
\begin{example}
\label{ex4}
Functions
$$
P=k_1 \left( k_2\frac{a_1-3/2}{a_2^3} -\frac{k_4}{a_2}+k_3 \right),
$$
$$
R=k_5 + \left( 4k_1k_4(2a_1-3)a_2^2 + k_1k_2 \left( 8a_0a_2-8a_2^2-3(3-2a_1)^2 \right) \right) / \left(8a_2^4 \right),
$$
$$
S=k_1k_4\log a_2-2k_5a_1+k_1k_3(a_2-2a_0)+k_6+ k_1\left(8a_2^4 \right)^{-1} \times
$$
$$
\times \left( k_2 \left( 2a_2^2(16a_1-27)+20a_0a_2(3-2a_1)+3(2a_1-3)^3 \right) +4k_4a_2^2 (6a_0a_2-4a_1(a_1-3)-9) \right)
$$
satisfy~\eqref{n3pde}, here $k_i\in\mathbb{R}$ are arbitrary constants. System~\eqref{n3eq} with such functions $P,$ $R,$ $S$ defines the implicit solution to~\eqref{n3system}.
\end{example}

\begin{remark}
\label{rem2}

In a typical situation integrable examples obtained via the generalized hodograph method turned out to be implicit (see Theorem~\ref{n3maintheorem} and Examples~\ref{ex1},~\ref{ex3},~\ref{ex4}). However, sometimes one can get rid of this implicitness via an appropriate change of coordinates (see Example~\ref{ex2} obtained from Example~\ref{ex1}). The same can be done for Examples~\ref{ex3},~\ref{ex4}. We skip these calculations.

\end{remark}

\section{Case $n=4$}

In this section we consider the case $n=4.$ We have $U=(a_0, a_1, a_2, a_3)^T,$ $a_4(t,x)\equiv 1,$ $g(t,x) = a_3(t,x),$
\begin{equation}
\label{genFn4}
F=\frac{a_0}{a_3^4}p_1^4+\frac{a_1}{a_3^3}p_1^3p_2+\frac{a_2}{a_3^2}p_1^2p_2^2+p_1p_2^3+p_2^4
\end{equation}
and system~\eqref{key_system_general} is of the form
\begin{equation}
\label{n4system}
U_t+VU_x=0, \qquad V=
\begin{pmatrix}
		0 & 0 & 0 & a_1 \\
		a_{3} & 0 & 0 & 2a_{2}-4a_{0} \\
		0 & a_{3} & 0 & 3a_3-3a_1 \\
    0 & 0 & a_{3} & 4-2a_2
\end{pmatrix}.
\end{equation}

As in the previous section, we apply the generalized hodograph method to system~\eqref{n4system} to construct its solutions. First of all we have to describe its symmetries. We have

\begin{lemma}
\label{n4lemma}
For matrix $V$ of the form~\eqref{n4system} the condition $[V,W]=0$ holds true if and only if $W$ has the form
\begin{equation}
\label{n4symm}
W=
\begin{pmatrix}
		w_1 & a_1P & a_1R & a_1S \\
		w_2 & w_3 & a_1P+2(a_2-2a_0)R & a_1R+2(a_2-2a_0)S \\
    a_3R+2(a_2-2)P & a_3S+2(a_2-2)R & T+2(a_2-2)S & w_4 \\
		a_3P & a_3R & a_3S & T
\end{pmatrix}
\end{equation}
with arbitrary $P,$ $R,$ $S,$ $T.$ Here
$$
w_1 = 2(2a_0-a_2)P+3(a_1-a_3)R+2(a_2-2)S+T,
$$
$$
w_2=3(a_1-a_3)P+2(a_2-2)R+a_3S, \qquad
w_3=3(a_1-a_3)R+2(a_2-2)S+T,
$$
$$
w_4=a_1P+2(a_2-2a_0)R+3(a_3-a_1)S.
$$

\end{lemma}

This lemma can be proved by direct calculations.

\begin{theorem}
\label{n4maintheorem}
Assume that functions $P(a_0,a_1,a_2,a_3),$ $R(a_0,a_1,a_2,a_3),$ $S(a_0,a_1,a_2,a_3)$ $T(a_0,a_1,a_2,a_3)$ satisfy the following system of PDEs:
\begin{equation}
\label{n4pde}
\begin{gathered}
a_1 \partial_0 P = 2(2a_0-a_2) \partial_0 R + 3(a_1-a_3) \partial_1 R + 2(a_2-2) \partial_2 R+a_3\partial_3 R+R,\hfill\\
\partial_1 P = \partial_0 R, \qquad \partial_2 P = \partial_1 R, \qquad \partial_3 P = \partial_2 R\hfill\\
\partial_0 S = \partial_1 R, \qquad \partial_1 S = \partial_2 R, \qquad \partial_2 S = \partial_3 R\hfill\\
a_3 \partial_3 S = a_1 \partial_0 R + 2(a_2-2a_0) \partial_1 R + 3(a_3-a_1) \partial_2 R-2(a_2-2)\partial_3 R+P,\hfill\\
\partial_0 T = a_3 \partial_2 R - 2P, \qquad \partial_1 T = a_3 \partial_3 R - 2R,\hfill\\
\partial_2 T = a_1 \partial_0 R + 2(a_2-2a_0) \partial_1 R + 3(a_3-a_1) \partial_2 R - 2(a_2-2)\partial_3 R + P-2S,\hfill\\
a_3 \partial_3 T = \gamma_0 \partial_0 R + \gamma_1 \partial_1 R + \gamma_2 \partial_2 R +\gamma_3 \partial_3 R + 2(2-a_2)P+2a_3R ,\hfill
\end{gathered}
\end{equation}
where $\partial_i = \frac{\partial}{\partial {a_i}},$ $i=0, 1, 2, 3$ and
$$
\gamma_0 = 2a_1(2-a_2), \qquad \gamma_1 = -4a_2^2+8a_0a_2+a_1a_3-16a_0+8a_2,
$$
$$
\gamma_2=-4a_0a_3+6a_1a_2-4a_2a_3-12a_1+12a_3, \qquad \gamma_3 = 4(a_2-2)^2-3a_1a_3+3a_3^2.
$$
Then any smooth solution $a_k(t,x)$ to the system of equations
\begin{equation}
\label{n4eq}
P=0, \qquad R=0, \qquad S=t, \qquad T=(4-2a_2)t-x
\end{equation}
is also a solution to system~\eqref{n4system}.

\end{theorem}

The proof of Theorem~\ref{n4maintheorem} is completely analogous to the one of Theorem~\ref{n3maintheorem}. Condition~\eqref{n4pde} is equivalent to~\eqref{matrixformsymmcond2} and relations~\eqref{n4eq} follow from system~\eqref{nondiaghodog} with $W$ of the form~\eqref{n4symm}. We skip the details.

Theorem~\ref{n4maintheorem} describes solutions to system~\eqref{n4system}. As in the case $n=3,$ it is difficult to construct a general solution to system~\eqref{n4pde}. So we shall restrict ourselves by constructing several partial solutions of a more or less simple form.

Assume that $P$ is a non-homogeneous quadratic polynomial in $a_0, \ldots, a_3$ with unknown constant coefficients. As previously, by substituting such $P$ into system~\eqref{n4pde} and integrating it, we find all the unknown functions $P(a_0, a_1, a_2, a_3), \ldots, T(a_0, a_1, a_2, a_3).$ This allows to obtain the following
\begin{example}
\label{ex5}
Functions $P,$ $R,$ $S,$ $T$ of the form
$$
5P= \left(64a_0^2+8a_1^2+4a_2^2+5a_3^2+16a_0a_2+8a_1a_3 \right)n_6-32a_0n_6+5(4a_0+a_2)n_2+5n_5,
$$
$$
R=\left(\frac{8}{5}(2a_0+a_2)n_6+n_2\right)a_1 + \left( \left( \frac{8}{5}a_0+2a_2+4 \right) n_6+3n_2/2 \right)a_3 + n_1/a_3,
$$
$$
S=\frac{4}{5}\left( 2a_0^2+a_1^2 + \frac{5}{4}a_2^2 + \frac{35}{8}a_3^2 +2a_0a_2+ \frac{5}{2}a_1a_3 \right)n_6+a_0n_2+4a_2n_6+\frac{3}{2}a_2n_2
$$
$$
+\left(6n_2+n_5+16n_6\right) \log a_3 + \frac{(2-a_2)}{a_3^2}n_1+n_3,
$$
$$
T = \frac{(2(a_2-2)^2-3a_1a_3)}{a_3^2}n_1 - \left( n_1+2(a_2-2)(6n_2+n_5+16n_6) \right) \log a_3 + \frac{\gamma}{60},
$$
where
$$
\gamma = 60n_4-64a_0(8a_0^2+3a_1^2-6a_0)n_6-40a_2^3n_6-6a_2^2(15n_2+8(2a_0+5)n_6)
$$
$$
-24a_2((8a_0^2+4a_1^2+5a_1a_3-40)n_6+5(a_0-3)n_2+5n_3)+60a_2n_5-60a_1^2n_2
$$
$$
-6a_1a_3(15n_2+8(2a_0+5)n_6) - 15 (16a_0^2n_2-3a_3^2(5n_2+32n_6)+8a_0n_5),
$$
satisfy~\eqref{n4pde}, here $n_i\in\mathbb{R}$ are arbitrary constants. System~\eqref{n4eq} with such functions $P,$ $R,$ $S,$ $T$ defines the implicit solution to~\eqref{n4system}.
\end{example}

Our next step is to extract various solutions of a simple form from Example~\ref{ex5}. Assuming that $n_5=-6n_2,$ $n_1=n_6=0$ and solving~\eqref{n4eq} w.r.t. $a_k(t,x)$ we obtain an explicit example which was constructed in~\cite{AbdMir}. We do not write it down here for brevity.

\begin{remark}
\label{rem3}

The integrable example in~\cite{AbdMir} was also obtained with the help of the generalized hodograph method but more straightforwardly. More precisely, in contrast with our approach, in~\cite{AbdMir} all the components of commuting flows are assumed to be polynomial in the field variables $a_k$ with unknown constant coefficients from the very beginning, and after that compatibility conditions~\eqref{nondiagcompcond} were verified straightforwardly. Recently the same "straightforward" approach yielded new integrable examples in slightly different situations (\cite{APSh},~\cite{AS2}).

We note that splitting the compatibility conditions~\eqref{nondiagcompcond} into two parts --- the algebraic one~\eqref{matrixformsymmcond1} and the differential one~\eqref{matrixformsymmcond2} --- essentially simplifies the search for symmetries and allows to extend the space of possible solutions.

\end{remark}

Let us go back to Example~\ref{ex5}. Assume now that $n_1=0,$ $n_2=k,$ $n_3= \ldots = n_6=0.$ Solving the first three equations of~\eqref{n4eq} w.r.t. $a_0(t,x)$, $a_1(t,x),$ $a_2(t,x)$ we are left with the last one equation on $a_3(t,x).$ To sum it up, we obtain

\begin{example}
\label{ex6}
Assume that $a_3(t,x)$ is a solution to the equation
$$
75k^2a_3^2+96k \left(6k\log a_3-2t-k \right) \log a_3+16t^2+16kt+20kx=0.
$$
Then the geodesic flow of the metric $ds^2 = g^2 (t,x) dt^2 + dx^2$ (recall that $a_3\equiv g$) admits the polynomial in momenta first integral $F_4$ of the form
$$
F_4=\frac{a_0}{a_3^4}p_1^4+\frac{a_1}{a_3^3}p_1^3p_2+\frac{a_2}{a_3^2}p_1^2p_2^2+p_1p_2^3+p_2^4,
$$
with coefficients
$$
a_0 = \frac{6k \log a_3-t}{5k}, \qquad a_1=-\frac{3a_3}{2}, \qquad a_2=-4a_0,
$$
here $k\in\mathbb{R}$ is an arbitrary constant.

\end{example}

An explicit example can be obtained from Example~\ref{ex5} in the following way. For simplicity assume that $n_6=0.$ Let us solve the first two equations of~\eqref{n4eq} $P=R=0$ w.r.t. $a_1(t,x),$ $a_2(t,x).$ Now, as usual, let us make the change of variables $(t,x) \rightarrow (a_0,a_3)$ and extend this transformation to the canonical one. Rewriting the metric and the first integral in terms of new coordinates and momenta we construct a family of explicit integrable examples. Assuming for simplicity $n_1=n_5=0,$ $n_2=1$ we obtain

\begin{example}
\label{ex7}
Denote $u^1=x,$ $u^2=y.$ The geodesic flow of the metric $ds^2 = g_{ij}(u)du^idu^j$ where
$$
g_{11} = 64(10x-1)^2+100y^2, \qquad g_{12} = 240y(5x-1), \qquad g_{22} = 9(25y^2+16)
$$
admits the fourth degree polynomial first integral $F=\alpha f_{k}(x,y)p_1^{4-k}p_2^k,$ where
$$
f_0 = 81(125(5x-6)y^4-320(5x-1)y^2+256),
$$
$$
f_1 = -27y(1875y^4+400(100x^2-80x+3)y^2-1024(50x^2-10x+3)),
$$
$$
f_2 = 72(125(80x-3)y^4+80(1500x^3-800x^2+95x+12)y^2-512(10x-1)^2x),
$$
$$
f_3= 24y(625y^4-200(700x^2-100x-7)y^2-128(10x-1)^2(100x^2-20x+3)),
$$
$$
f_4 = 4096(10x-1)^4x+2560(10x-1)^2(20x-3)y^2-2000(40x-9)y^4,
$$
here $\alpha = (25y^2+160x-16)^{-4}.$ Gauss curvature $K$ of this metric is equal to
$$
K = \frac{25 (25y^2-160x+208)}{9 (25y^2+160x-16)^3}.
$$

The geodesic flow of this metric does not admit linear in momenta first integrals.

\end{example}

As in the previous case $n=3,$ many other integrable examples with an integral of degree $n=4$ can be constructed by choosing $P$ in a polynomial form of higher degrees.

\section{Case $n=5$}

In this section we consider the case $n=5.$ We follow the same algorithm as in the cases $n=3,4$ but calculations and the final form of equations and functions become very complicated. Therefore we shall briefly describe only the main points.

We have $U=(a_0, a_1, a_2, a_3,a_4)^T,$ $a_5(t,x)\equiv 1,$ $g(t,x) = a_4(t,x),$
\begin{equation}
\label{genFn5}
F=\frac{a_0}{a_4^5}p_1^5+\frac{a_1}{a_4^4}p_1^4p^2+\frac{a_2}{a_4^3}p_1^3p_2^2+\frac{a_3}{a_4^2}p_1^2p_2^3+p_1p_2^4+p_2^5
\end{equation}
and system~\eqref{key_system_general} is of the form
\begin{equation}
\label{n5system}
U_t+VU_x=0, \qquad V=
\begin{pmatrix}
		0 & 0 & 0 & 0 & a_1 \\
		a_{4} & 0 & 0 & 0 & 2a_{2}-5a_{0} \\
		0 & a_{4} & 0 & 0 & 3a_3-4a_1 \\
    0 & 0 & a_{4} & 0 & 4a_4-3a_2\\
    0 & 0 & 0 & a_{4} & 5-2a_3 \\
\end{pmatrix}.
\end{equation}

\begin{lemma}
\label{n5lemma}
For matrix $V$ of the form~\eqref{n5system} the condition $[V,W]=0$ holds true if and only if $W$ has the form
$$
W=
\begin{pmatrix}
		* & a_1P & a_1R & a_1S & a_1T \\
		* & * & * & * & * \\
    * & * & * & * & * \\
		* & * & * & * & * \\
    a_4P & a_4R & a_4S & a_4T & Q
\end{pmatrix}
$$
with arbitrary $P,$ $R,$ $S,$ $T,$ $Q.$ All the missing components denoted by $*$ are of a very cumbersome form, we skip them for brevity.
\end{lemma}

\begin{theorem}
\label{n5maintheorem}
Assume that functions $P(a_0,a_1,a_2,a_3,a_4), \ldots, Q(a_0,a_1,a_2,a_3,a_4)$ satisfy a system of PDEs similar to the ones~\eqref{n3pde},~\eqref{n4pde} in cases $n=3, 4$ (we skip its explicit form for brevity). Then any smooth solution $a_k(t,x)$ to the system
$$
P=0, \qquad R=0, \qquad S=0, \qquad T=t, \qquad Q=(5-2a_3)t-x
$$
is also a solution to system~\eqref{n5system}.

\end{theorem}

As previously, various integrable examples can be constructed by choosing $P$ in a polynomial form. Assume that $P$ is linear in $a_0, \ldots, a_4.$ Finding $P, \ldots, Q$ from condition~\eqref{matrixformsymmcond2} we obtain a family of integrable examples parameterized by a number of arbitrary constants. The simplest example has the following form.

\begin{example}
\label{ex8}
Assume that functions $a_0(t,x),$ $a_1(t,x),$ $a_2(t,x)$ $a_3(t,x),$ $a_4(t,x)$ satisfy
\begin{equation}
\label{systex8}
\begin{gathered}
24a_0+4a_2+3a_4=0, \qquad 8a_1+6a_3+15=0,\hfill\\
8a_0+6a_2+15a_4=0, \qquad 6a_1+15a_3+105 \log a_4=8t/k,\hfill\\
96a_0^2+16a_1^2+12a_2^2-30a_3^2-105a_4^2+32a_0a_2+12a_0a_4+30a_2a_4+120a_1-60a_3=16x/k,\hfill
\end{gathered}
\end{equation}
here $k\in\mathbb{R}$ is an arbitrary constant.

Then the geodesic flow of the metric $ds^2 = g^2 (t,x) dt^2 + dx^2$ (recall that $a_4\equiv g$) admits the fifth degree polynomial in momenta first integral $F$ of the form~\eqref{genFn5}.
\end{example}

In order to rewrite this example in an explicit form let us solve the first three equations of~\eqref{systex8} w.r.t. $a_0,$ $a_1,$ $a_2.$ Then, as previously, let us make the change of coordinates $(t,x) \rightarrow (a_3, a_4)$ via the last two equations of~\eqref{systex8} and extend it to the canonical transformation. Rewriting the metric and the first integral in terms of new coordinates and momenta we obtain

\begin{example}
\label{ex9}
Denote $u^1=x,$ $u^2=y.$ The geodesic flow of the metric $ds^2 = g_{ij}(u)du^idu^j$ where
$$
g_{11} = (2x+5)^2+y^2, \qquad g_{12} = 20y(x+3), \qquad g_{22} = 100(y^2+1)
$$
is completely integrable. Gauss curvature $K$ of this metric is equal to
$$
K = \frac{(y^2+2x+25)}{100(y^2-2x-5)^3}.
$$
The additional fifth degree polynomial first integral is the following:
$$
F=\frac{1}{(y^2-2x-5)^5} \left( 3y\alpha^5-3(2x+5)\alpha^4\beta-24y\alpha^3\beta^2+8x\alpha^2\beta^3+8y\alpha\beta^4+8\beta^5 \right),
$$
where $\alpha = 10yp_1-(2x+5)p_2,$ $\beta = 10p_1-yp_2.$

This geodesic flow does not admit any linear integrals.
\end{example}

\section{Conclusion}

This paper deals with the classical problem of searching for 2-dimensional Riemannian metrics with polynomially integrable geodesic flows. As shown in~\cite{BM1}, this problem reduces to searching for solutions to semi-Hamiltonian system of PDEs~\eqref{key_system_general},~\eqref{key_matrix_general}. The most difficult part of the present work is the one related to searching for commuting flows of~\eqref{key_system_general},~\eqref{key_matrix_general}. We write down the full set of PDEs on the components of commuting flows (see systems~\eqref{n3pde},~\eqref{n4pde} for $n=3, 4$ correspondingly) and construct plenty of its partial solutions. Finally, applying the generalized hodograph method we construct a series of explicit and implicit examples of local 2-dimensional metrics with integrable geodesic flows admitting additional polynomial in momenta first integrals of degrees 3, 4 or 5 (see Examples~\ref{ex1} ---~\ref{ex9}).

Let us make some concluding remarks.

1) We believe that following this approach, it is possible to construct new integrable examples of 2-dimensional geodesic flows with polynomial integrals of degrees $n>5.$ However, all the necessary calculations and the final form of metrics and first integrals in this case are expected to be rather cumbersome.

2) We notice that, despite of complicated calculations, in all the cases $n=3, 4, 5$ the system of equations~\eqref{nondiaghodog} turns out to have a very simple form
$$
S_1 = S_2 = \ldots = S_{n-2} = 0, \qquad S_{n-1} = t, \qquad S_n = (n-2a_{n-2})t - x,
$$
where $n$ functions $S_k(a_0, \ldots, a_{n-1})$ determine all the components of commuting flows. Moreover, the structure of these equations will be the same for any $n \in \mathbb{N}.$

3) It is a very interesting problem to construct a general solution to systems~\eqref{n3pde},~\eqref{n4pde}. It would allow to describe a wider class of integrable metrics admitting polynomial integrals of degrees 3, 4.

4) In a typical situation examples obtained via the generalized hodograph method turn out to be local. It would be very interesting to understand whether it is possible to construct smooth periodic solutions on the 2-torus or to prove non-existence of such solutions (see Conjecture~\ref{Conj} in Introduction) following this approach.

5) In this paper we work with semi-geodesic coordinates $ds^2=g^2(t,x)dt^2+dx^2.$ After applying the generalized hodograph method we obtain solutions in an implicit form, i.e. $g(t,x)$ and the rest of unknown functions satisfy certain systems of algebraic or transcendental equations (e.g., see~\eqref{ex1impl1},~\eqref{ex3impl1} and others). In principle, it may happen that for some specific values of free parameters function $g$ turns out to be purely imaginary, which gives rise to Lorentzian metrics. We did not investigate this question in the paper.

{\bf Acknowledgements.} The author thanks S.P. Tsarev and M.V. Pavlov for helpful discussions. The author also thanks anonymous reviewers for their comments and suggestions which helped to significantly improve the text.

${\ }$

Sergei Agapov,

Novosibirsk State University, 630090, 1, Pirogova str., Novosibirsk;

Sobolev Institute of Mathematics SB RAS,  630090, 4 Acad. Koptyug avenue, Novosibirsk

e-mail: agapov.sergey.v@gmail.com

\end{document}